\documentclass{amsart}
\usepackage{amsfonts}
\usepackage{amssymb}
\usepackage{amscd}

\newcommand{\Z}{\mathbb{Z}}
\newcommand{\dbR}{\mathbb{R}}
\newcommand{\dbK}{\mathbb{K}}
\newcommand{\dbE}{\mathbb{E}}

\newcommand{\dbH}{\mathbb{H}}

\newcommand{\ca}{{\mathcal A}}
\newcommand{\cb}{{\mathcal B}}
\newcommand{\cf}{{\mathcal F}}
\newcommand{\cg}{{\mathcal G}}
\newcommand{\bfK}{{\mathbf K}}

\newcommand{\borG}{{\Omega_n^{\rm Spin}(BG)}}

\newcommand{\koG}{{ko_n(BG)}}

\newcommand{\KoG}{{KO_n(BG)}}

\newcommand{\cstarc}{{C^*\mbox{-CATEGORIES}}}

\newcommand{\cstarG}{{KO_n(C_r^*G)}}

\newcommand{\G}{{\Gamma}}

\newcommand{\fin}{{\mathcal F}in}

\newcommand{\Parrow}{\buildrel p \over \longrightarrow}
\newcommand{\Aarrow}{\buildrel A \over \longrightarrow}
\newcommand{\Darrow}{\buildrel D \over \longrightarrow}

\DeclareMathOperator{\id}{Id}

\DeclareMathOperator{\topo}{top}
\DeclareMathOperator{\proper}{proper}
\DeclareMathOperator{\Or}{Or}

\DeclareMathOperator{\groupoids}{GROUPOIDS}
\DeclareMathOperator{\spectra}{SPECTRA}
\DeclareMathOperator{\sets}{SETS}
\DeclareMathOperator*{\hocolim}{hocolim}
\DeclareMathOperator*{\mor}{mor}
\newtheorem{theorem}{Theorem}

\newtheorem{proposition}[theorem]{Proposition}
\newtheorem{corollary}[theorem]{Corollary}

\newtheorem{remark}[theorem]{Remark}

\newtheorem*{main}{Main Theorem}
\newtheorem*{GLR}{The Gromov-Lawson-Rosenberg (GLR) Conjecture for
$G$}
\newtheorem*{bordism}{Bordism Theorem}
\newtheorem*{obstruction}{Obstruction Theorem}

\begin{document}
\title[The GLR conjecture for cocompact Fuchsian groups]{The
Gromov-Lawson-Rosenberg conjecture for cocompact Fuchsian groups}

\author{James F. Davis}
\address{Department of Mathematics\\
Indiana University\\
Bloomington, IN 47405}
\email{jfdavis@indiana.edu}

\author{Kimberly Pearson}
\address{Department of Mathematics and Computer Science\\
Valparaiso University\\
Valparaiso, IN 46383}
\email{kpearson@valpo.edu}

\thanks{{\it 2000 Mathematics Subject Classification} Primary: 53C21.
Secondary: 19L41, 19L64, 57R15, 55N15, 53C20.}
\thanks{{\it Key words and phrases:} positive scalar curvature, $K$-theory,
Fuchsian groups}

\begin{abstract}
We prove the Gromov-Lawson-Rosenberg conjecture for cocompact Fuchsian
groups, thereby giving necessary and sufficient conditions for a closed
spin manifold of dimension greater than four with
fundamental group cocompact Fuchsian to admit a metric of positive
scalar curvature.
\end{abstract}

\maketitle

Given a smooth closed manifold $M^n$, it is a long-standing question to
determine
whether or not $M$ admits a Riemannian metric of positive
scalar curvature.
Work of Gromov-Lawson and
Schoen-Yau shows that if $N^n$ admits positive scalar curvature
and $M$ is obtained
from $N$ by $k$-surgeries of codimension $n-k\geq 3$,
then $M$ admits positive scalar curvature as well.
In the case when $M$ is spin, this surgery result implies
the following.
\begin{bordism} (\cite{GL}, \cite{SY})  Let $M^n$ be a closed spin
manifold, $n \geq
5$,
 $G = \pi_1(M)$, and suppose $u:M\to BG$ induces the  identity on the
fundamental group.  If there is a positively  scalar
curved spin
manifold $N^n$
 and a map $v: N \to BG$
such that $[M,u] = [N,v] \in \borG$, then $M$ admits a metric of positive
scalar curvature.
\end{bordism}
On the other hand, the work of Lichnerowicz gives an obstruction
to manifolds admitting positive scalar curvature.
Using the Weitzenb\"ock formula for the Dirac operator and the
Atiyah-Singer index theorem, he
proves in \cite{Li} that if $M^{4k}$ is a closed spin manifold with
positive scalar curvature then $\hat A(M)$, the $A$-hat genus of
$M$, vanishes.
Generalizations of the Dirac operator and its index by
Hitchin \cite{Hi} and Mi{\v{s}}{\v{c}}enko-Fomenko \cite{MF}
provide obstructions as well,
taking final form in the following theorem of  Rosenberg \cite{R3}.
\begin{obstruction}\label{Rose}\cite{R3} \cite{S3} Let $M^n$ be a closed spin
manifold and
$u:M\to BG$ be a continuous map for some discrete group $G$.   If $M$
admits a metric of
positive scalar curvature then $\alpha[M,u] = 0$ in $\cstarG$, where $\alpha
: \borG \to \cstarG$ is the index of the Dirac
operator.
\end{obstruction}

Here $C^*_rG$ is the reduced real $C^*$-algebra of $G$,
a suitable completion of the group algebra $\dbR G$.
In dimensions $4k$ with $G$ = 1, the index  $\alpha[M,u]$ agrees with $\hat
A(M)$
up to a constant factor (\cite{LM}, p.~149).  The above theorem
motivates
\begin{GLR} Let $M^n$ be a\linebreak closed spin manifold,
$n\geq 5$, $G = \pi_1(M)$,
and suppose $u:M\to BG$ induces the  identity on the fundamental group.
Then $M$ admits positive scalar curvature  if
and only if
$\alpha[M,u]$ is zero in $KO_n(C^*_rG)$.
\end{GLR}
\noindent
The conjecture has been verified in the cases when $G$ is
trivial \cite{S2}, $\Z/2$ \cite{RS}, $\Z/n$ for $n$ odd (\cite{KS}, \cite{R2}),
has periodic cohomology (\cite{BGS}, \cite{KS}), or is one of certain
elementary abelian groups \cite{BG}.  Infinite groups for which
the conjecture is proven are free and free abelian groups \cite{RS2},
certain other torsion-free groups \cite{JS}, and
some crystallographic groups with prime torsion \cite{DL3}.
However, in \cite{Sch},
Schick shows the conjecture is false for $\G = \Z^4\oplus \Z/3$, and other
counterexamples have since been
constructed in \cite{JS}.  In \cite{RS2},
Rosenberg and Stolz introduced a ``stable'' version:
that the conjecture holds after crossing $M$ with sufficiently
many copies of the Bott manifold $B^8$.

In this paper we prove the original (unstable) conjecture when
$G$ is a cocompact Fuchsian group $\G$.  Fuchsian groups are the discrete
subgroups
of $PSL_2(\dbR)$, which can be identified with the orientation-preserving
isometries of the hyperbolic plane $\dbH^2$, and as such they
have rich geometric properties (see \cite{K} as a reference).
Each cocompact Fuchsian group is classified by its signature, $(g; \nu_1,
\nu_2, ..., \nu_r)$, where
$g$ is the genus of the compact orbit space $\dbH^2/\G$ and $\Z/\nu_j$ is
the stabilizer of a lift $\tilde v_j$ of a singular point $v_j\in \dbH^2/\G$ of
the branched covering $\dbH^2 \to \dbH^2/\G$.
The  group
 $\G$ has $r$ conjugacy classes of maximal finite subgroups, all of which are
cyclic. The integers $\nu_i$ indicate their order.
Furthermore, nearly all signatures $(g; \nu_1, \nu_2,\ldots, \nu_r)$ are
realizable: if
$\nu_j\geq 2, g\geq 0, r\geq 0$ and $2g -2 + \Sigma (1 - 1/\nu_j) > 0$,
then there exists a Fuchsian group $\G$ with this signature. Generators and
relations for $\G$ were found by Poincar\'e \cite{P} (a more modern
reference is
\cite{W}),
$$
\G =\langle a_1, b_1, \dots a_g, b_g, c_1, \dots c_r~|~ [a_1,b_1]\cdots
[a_g,b_g] = c_1 \cdots c_r,
c_1^{\nu_1} = \cdots = c^{\nu_r}= 1  \rangle .
$$
We can now formally state our
\begin{main}
The Gromov-Lawson-Rosenberg conjecture is true for all cocompact Fuchsian
groups.
\end{main}
Following the general line of proof in \cite{DL3}, our strategy is to use the
validity of Baum-Connes
conjecture for the group $\G$ and  the $p$-chain spectral sequence
to
understand the map $\alpha$ by
studying how it factors through connective and
periodic real $K$-theory, and then to exploit the fact that the
Gromov-Lawson-Rosenberg conjecture is true for each
finite subgroup of $\G$.

\section{Assorted Definitions}
A real $C^*$-algebra is a real Banach $*$-algebra $\mathcal A$ which
is $*$-isomorphic to a norm-closed subalgebra of the bounded
operators $\cb(H)$ on a real Hilbert space $H$.  An example of such is
given by the reduced $C^*$-algebra $C^*_rG$ of a discrete
group
$G$, which is defined to the closure of $\dbR G$ in $\cb(\ell^2(G))$,
where  $\dbR G$ acts on $\ell^2(G)$ by multiplication.

For a $C^*$-algebra $\ca$ with unit, and for $n > 0$, define
$KO_n(\ca) = \pi_{n-1}(GL(\ca))$.  There is the Bott map, defined
using Clifford algebras, $GL(\ca) \to \Omega^8(GL(\ca))$, which,
according to Bott periodicity \cite{Wood}, is a homotopy
equivalence.  Then for any $n$, define $KO_n(\ca) = \pi_{n + 8k-1}(GL(\ca))$,
where $n + 8k > 0$.  One can show
 $KO_0(\ca)$ is the Grothendieck group of isomorphism
classes of finitely generated
projective $\ca$-modules, in other words, the algebraic
$K_0$ of the ring $\ca$.  For more details on $C^*$-algebras
and their $K$-theory we refer the reader to \cite{LM}.

The above theory is sufficient for analysis, but to compute
using algebraic topology one needs spectra.  By a spectrum
$\bfK$, we mean an sequence of based spaces
$\{K_n\}_{n\geq 0}$ and based maps
$\sigma_n :S K_n
\to K_{n+1}$, where $S$ denotes suspension.  A map of spectra
$f : \bfK \to \bfK'$ is a
sequence of based maps $f_n : K_n \to K_n'$ commuting with the structure
maps, i.e.
$f'_{n+1} \circ \sigma_n = \sigma'_n \circ Sf_n$.  For a real
$C^*$-algebra with unit $\ca$, define the Bott spectrum ${\bf
KO(\ca)}$ to have as its $n$-th space $\Omega^i(GL(\ca))$
where $i \in \{0,1,2,\dots, 7\}$ and $i+n +1  \equiv 0 \pmod 8$.
Most
of the structure maps $S\Omega^i(GL(\ca)) \to \Omega^{i-1}(GL(\ca))$
are given as the adjoint of the identity map on
$\Omega^i(GL(\ca))$, while one-eighth of the maps are given by the
adjoint of the Bott map.  Then $\pi_n ({\bf KO(\ca)}) =
KO_n(\ca)$ for all $n \in \Z$.  We will write ${\bf KO}$ for
${\bf KO}(\dbR)$.

However, the above spectra are still insufficient for our purposes.
A preliminary problem is that the fundamental group depends
on a choice of base point, but a much more severe problem is that
$(X,x_0) \mapsto C^*_r(\pi_1(X,x_0))$ is not a
functor from the category of based spaces, since the reduced
$C^*$-algebra of a group is not functorial for group
homomorphisms with infinite kernel.  This obstructs both a
definition of the assembly map in terms of algebraic topology
and an axiomatic characterization of the assembly map, both
of which are necessary in this paper.

These base point and functorial problems can be avoided by the
use of $C^*$-categories, which is the point of view of
\cite{DL1}.  In this paper we will need a composite of three
functors
$$
\Or(G) \to \groupoids^{\proper} \to \cstarc \to \spectra
$$
which we describe from right to left.  A $C^*$-category $\ca$ is a
category, so that for any
two objects $c,d$, the morphism set  $\mor_\ca(c,d)$ is a Banach
space with an involution $* : \mor_\ca(c,d) \to
\mor_\ca(d,c)$ satisfying the expected properties (e.g.
bilinearity with respect to composition), and so that for all $f
\in
\mor_\ca(c,d)$, $\|f^*
\circ f\| =
\|f\|^2$ and so that $f^* \circ f$ has non-negative spectrum
(in the sense of analysis).  A $C^*$-category is a topological
category: the morphisms sets are topological spaces, although
in this case the objects are discrete.  A
$C^*$-category with exactly one object is a $C^*$-algebra. The
category
$\cstarc$ is the category of all small (i.e. the objects form a set)
$C^*$-categories.  The functor $F : \cstarc \to \spectra$ is
to have two properties: when the $C^*$-category has a single object,
one gets the Bott spectrum and and an equivalence of
$C^*$-categories gives a  homotopy equivalence of spectra.
Unfortunately, the construction of $F$ in \cite{DL1} was
flawed by failing to heed the warning of \cite{T}; the problem
can be repaired without a great deal of difficulty using Segal's
$\Gamma$-spaces.  Alternatively, such a functor has been constructed in
\cite{M}.

A (discrete) groupoid $\cg$ is a category all of whose morphisms
are isomorphisms; a groupoid with a single object is a group.
Associated to a groupoid $\cg$ is a $C^*$-category $C^*_r(\cg)$, whose
objects are the same as those of $\cg$, and whose morphisms
are given by the closure of the linear span of
$\mor_\cg(c,d)\subset \cb(\ell^2(\mor_\cg (c,c), \ell^2(\mor_\cg
(c,d))$.  A morphism of groupoids is a functor; a
morphism $F : \cg_1 \to \cg_2$
is proper if for all objects $c$ and $d$
of
$\cg_1$, there is a $n = n(c,d)$ so that for all
$f \in \mor_{\cg_2}(F(c),F(d))$, $F^{-1}(f) \subset \mor_{\cg_1}(c,d)$
has cardinality less than $n$.  A proper morphism between groupoids induces a
map of the associated $C^*$-categories.  This gives our
functor from\linebreak
$\groupoids^{\proper}$, the category of small groupoids and proper
morphisms, to $\cstarc$.

Finally, for a group $G$, define the orbit category $\Or(G)$ whose
objects are the $G$-sets $G/H$ where $H$ is a subgroup of $G$
and whose morphisms are the equivariant maps.  The orbit category is a
useful indexing set in group actions, for
example, if $X$ is a $G$-set there is an associated functor
$\Or(G) \to \sets$, $G/H \mapsto X/H$.  For any $G$-set $X$,
there is an associated groupoid $\overline X$ whose objects
are the elements of $X$ and $\mor_{\overline
X}(x_1,x_2) = \{g \in G : gx_1 = x_2\}$, with the composition
law given by group multiplication.  This gives our last
functor
$\Or(G)
\to
\groupoids^{\proper}$.  Note the
$\overline {G/H}$ is equivalent as a category to the group $H$.

The composite functor
$$
{\bf KO}^{\topo}: \Or(G) \to \spectra,
$$
satisfies $\pi_n {\bf KO}^{\topo}(G/H) = KO_n(C^*_rH)$.
This functor is, according to \cite{DL1}, the basic building block
of the assembly map (of which there are four equivalent
descriptions in \cite{DL1}).  Give a family
$\cf$ of subgroups of $G$, define the restricted orbit category
$\Or(G,\cf)$ to be the full subcategory of $\Or(G)$ whose objects
are $G/H$ with $H \in \cf$.  Let $1$ denote the family consisting
of the trivial subgroup and $\fin$ the family of finite
subgroups.  Consider the maps
$$
BG_+ \wedge {\bf KO} \simeq\hocolim_{\Or(G,1)} {\bf KO}^{\topo}
\to \hocolim_{\Or(G,\fin)} {\bf KO}^{\topo}
\to
\hocolim_{\Or(G)} {\bf KO}^{\topo} \simeq {\bf KO}(C^*_rG)
$$
The assembly map $KO_n(BG) \to KO_n(C^*_r(G))$ is, by definition,
$\pi_n$ applied to the composite.  The
induced map
$$\pi_n(\hocolim_{\Or(G,\fin)} {\bf KO}^{\topo})
\to KO_n(C^*_r(G))$$
has been identified with the  Baum-Connes map
in \cite{HP}.  This identification is crucial for us, because
the analytic map gives the obstructions to positive scalar
curvature, and the topological map we need for computations.

\section{Proof of the Main Theorem}

Recall \cite{S3} that   $\alpha$ is the composition
$\alpha = A\circ p\circ D$,
\[
\borG\Darrow\koG\Parrow\KoG\Aarrow\cstarG .
\]
Here $D$ is the $ko$-orientation of spin bordism,
$p$ is the covering of periodic $K$-theory by  connective $K$-theory, and
$A$ is the
assembly map.  The groups $ko_n(*) = 0$ for $n < 0$ and
$p : ko_n(*) \to KO_n(*)$ is an isomorphism for
$n \geq 0$.

The utility of this factorization is the generalization
of the Bordism Theorem
by Jung and Stolz.  Let $ko_n^+(BG)$ be the
subgroup of $ko_n(BG)$ given by $D[N^n,v]$ where $N$ is a positively curved
spin manifold and $v : N \to BG$ is a continuous map.

\begin{theorem}  \cite{S3}  \label{kobordism} Let $M^n$ be a closed
spin manifold, $n \geq 5$,
 $G = \pi_1(M)$, and suppose $u:M\to BG$ induces the  identity on the
fundamental group.  If $D[M,u] \in ko_n^+(BG)$, then $M$
admits a metric of positive  scalar curvature.

\end{theorem}

  Our direct sum in the theorem below and hereafter
will be over the $r$ conjugacy
classes $(C)$ of maximal finite subgroups of $\G$.  Thus  $C \cong
\Z/\nu_j$ for  $j = 1, \dots r$.

\begin{proposition}\label{pcss}\cite{BJP} Let $\G$ be a cocompact Fuchsian
group with
signature\break $(g; \nu_1, \nu_2,\ldots, \nu_r)$
and $\dbE$ be any covariant functor from {\rm Or}$(\G,\fin)$ to {\rm
SPECTRA}.
Then there is an exact sequence
\[
\dots \rightarrow\oplus_{(C)} H_n(BC;\dbE(\G/1))\to (\oplus_{(C)}
\pi_n(\dbE(\G/C))
\oplus H_n(B\G;\dbE(\G/1))
\]
\[
\ \ \ \ \to\pi_n (\underset{{\rm Or}(\G,\fin)}{\rm hocolim}(\dbE ))\to
\oplus_{(C)} H_{n-1}(BC;\dbE(\G/1))\to \cdots.
\]
\end{proposition}
\noindent

This result is obtained by showing the spectral sequence
from \cite{DL2} collapses at $E^2$.  The proposition in \cite{BJP} is stated
for the algebraic $K$-theory functor $\dbK$, but holds for any $\dbE$. The
reason why
cocompact Fuchsian groups  are so well-behaved is due to the structure of
the lattice of
finite subgroups:  if $H$ is a non-trivial finite subgroup, then $H$ is
contained
in a  unique maximal finite subgroup
$M$, the normalizer $N_\G(H)$ equals $M$, and the set of finite subgroups
containing $H$ is totally ordered.

When $\dbE$ is a spectrum, and $\dbE_c:{\rm Or}(\G,\fin)\to \spectra$ is
the  constant functor
$\G/H \mapsto \dbE$,
the following simplification to Proposition \ref{pcss} can be made.

\begin{corollary}\label{cor}
Let $\G$ be a cocompact Fuchsian group with signature
$(g; \nu_1, \nu_2,\ldots, \nu_r)$, $\dbE$ be a spectrum, and
$X_g$ be the closed orientable surface of genus $g$.
Then
\[
\cdots \rightarrow\oplus_{(C)} H_n(BC;\dbE)\to (\oplus_{(C)}
\pi_n(\dbE))\oplus H_n(B\G;\dbE)
\to H_n(X_g;\dbE) \to \cdots
\]
and
\[
\cdots \rightarrow\oplus_{(C)} \widetilde H_n(BC;\dbE)\to  \widetilde
H_n(B\G;\dbE)
\to\widetilde H_n(X_g;\dbE) \to \cdots
\]
are exact.
\end{corollary}

\begin{proof}
By Proposition \ref{pcss}, it suffices to identify
$\pi_n (\hocolim_{\Or(\G,\fin)}(\dbE_c ))$ with\break $H_n(X_g;\dbE)$.
Let $E_{\fin}\G$ be the ``universal $\G$-space with finite isotropy'', i.e. a
$\G$-complex
such that the fixed set $(E_{\fin}\G)^H$ is contractible if $H$ is a finite
subgroup of $\G$
and empty if $H$ is not finite (see \cite{E}, \cite{tD}).
By definition of the homotopy colimit (see \cite[sections 3 and 7]{DL1}),
$\pi_n (\hocolim_{\Or(\G,\fin)}(\dbE_c ))\cong
H_n(E_{\fin}\G/\G;\dbE)$ for any constant Or($\G,\fin$)-spectrum
$\dbE_c$.
A cocompact Fuchsian group acts with
finite isotropy on the hyperbolic plane $\dbH^2$ with quotient $X_g$,
and we can take $\dbH^2$ as a model for $E_{\fin}\G$.
\end{proof}

\begin{proposition} \label{ladder} Let $\G$ be a cocompact Fuchsian group
with signature\linebreak
$(g; \nu_1, \nu_2,\ldots, \nu_r)$ and
$X_g$ be the closed orientable surface of genus $g$.  There is a commutative
diagram with exact rows
$$
\begin{CD}  \widetilde{KO}_{n+1}(X_g) @>>> \oplus_{(C)}\widetilde{KO}_n(BC)
@>>>  \widetilde{KO}_n(B\G) @>>>
\widetilde{KO}_{n}(X_g)  \\
@VV\id V @VV A  V @VV A V @VV\id V \\
\widetilde{KO}_{n+1}(X_g) @>>> \oplus_{(C)}\widetilde{KO}_n(C^*_rC) @>>>
\widetilde{KO}_n(C^*_r\G) @>>>
\widetilde{KO}_{n}(X_g).
\end{CD}
$$
\end{proposition}

\begin{proof}  In \cite{HP}, Hambleton-Pedersen identify
the Baum-Connes map with the Davis-L\"uck assembly map \cite{DL1}
$$
\pi_n\hocolim_{\Or(G,\fin)} {\bf KO}^{\topo} \to KO_n(C^*_rG).
$$
Kasparov \cite{Kas} \cite[p. 284]{BCH} proved
that Baum-Connes map is an isomorphism for subgroups of
$SO(n,1)$, such as $\G$.  Thus we can (and will) replace
$KO_n(C^*_r\Gamma)$ by $\pi_n\hocolim_{\Or(G,\fin)} {\bf KO}^{\topo}$ in the
proof of this proposition.

Underlying Proposition \ref{pcss} is an ``excision'' isomorphism,
$$
\pi_*(\vee_{(C)}BC_+ \wedge \dbE(\Gamma/1) \to B\Gamma_+ \wedge
\dbE(\Gamma/1) ) \xrightarrow{\simeq}
\pi_*(\vee_{(C)}\dbE(\Gamma/C) \to \hocolim_{\Or(\Gamma,\fin)}\dbE),
$$
defined for every functor $\dbE : \Or(\Gamma,\fin) \to \spectra  $.
Now notice that that for $H \in \fin$, there is a morphism in
$\groupoids^{\proper}$
$$
\overline{\Gamma/H} \to \overline{1/1},
$$
where the target is the trivial groupoid.  When $H = 1$, this is
an equivalence.  Hence there is a
natural transformation of $\Or(\Gamma,\fin)-\spectra$,
$$
{\bf KO}^{\topo} \to {\bf KO}_c,
$$
inducing a homotopy equivalence ${\bf KO}^{\topo}(\Gamma/1) \to
{\bf KO}_c(\Gamma/1) = {\bf KO}$.

  Proposition
\ref{ladder} now follows from a diagram chase.  Indeed,
consider the following diagram.
$$\begin{CD}
\pi_*(\vee_{(C)}BC_+ \wedge {\bf KO}\to B\Gamma_+ \wedge
{\bf KO} ) @>\simeq>>
\pi_*({\bf KO} \to X_{g+} \wedge {\bf KO}) \\
@VV\simeq V @VV\id V\\
\pi_*(\vee_{(C)}{\bf KO}(C^*_rC) \to {\bf KO}(C^*_r\Gamma) ) @>>>
\pi_*({\bf KO} \to X_{g+} \wedge {\bf KO})
\end{CD}$$
The top row is the excision isomorphism coming from ${\bf KO}_c$
and the left column is the excision isomorphism coming from ${\bf KO^{\topo}}$.
Commutativity of the diagram comes from the natural transformation
${\bf KO}^{\topo} \to {\bf KO}_c$.  The horizontal  isomorphisms
give Mayer-Vietoris exact sequences and the commutativity of the square gives
a map between the two Mayer-Vietoris sequences.
\end{proof}

\begin{proof}[Proof of Main Theorem]  Let ${\bf ko}$ be the
spectrum for the connective cover of ${\bf KO}$.  Since there is a
natural transformation $p : {\bf ko}_c \to {\bf KO}_c$ of
$\Or(\Gamma,\fin)$-$\spectra$, Corollary  \ref{cor} and Proposition
\ref{ladder} give a commutative diagram with exact rows
$$
\begin{CD}  \widetilde{ko}_{n+1}(X_g) @>>> \oplus_{(C)}\widetilde{ko}_n(BC)
@>>>  \widetilde{ko}_n(B\G) @>>>
\widetilde{ko}_{n}(X_g)  \\
@VVp V @VV A \circ p  V @VV A \circ p V @VVp V \\
\widetilde{KO}_{n+1}(X_g) @>>> \oplus_{(C)}\widetilde{KO}_n(C^*_rC) @>>>
\widetilde{KO}_n(C^*_r\G) @>>>
\widetilde{KO}_{n}(X_g).
\end{CD}
$$

Now $\Sigma X_g \simeq (\vee_{2g}~S^2) \vee S^3$, so
\begin{align*}
\widetilde{ko}_n(X_g) &\cong  ko_{n-1}(*)^{2g} \oplus ko_{n-2}(*) \\
\widetilde{KO}_n(X_g) &\cong  KO_{n-1}(*)^{2g} \oplus KO_{n-2}(*),
\end{align*}
and hence $p : \widetilde{ko}_n(X_g) \to \widetilde{KO}_n(X_g)$ is an
isomorphism for $n \geq 2$; an alternative proof uses the
Atiyah-Hirzebruch spectral sequence.
Suppose now that $n \geq 5$ and that
$$
\beta \in \ker (A \circ p : ko_n(B\G) \to KO_n(C^*_r\G)) \cong\ker (A \circ p :
\widetilde{ko}_n(B\G) \to \widetilde{KO}_n(C^*_r\G)).
$$
Then using the commutative diagram and the fact that the $p$'s are
isomorphisms, one finds an element
$$
\gamma \in \ker (A \circ p : \oplus_{(C)} ko_n(BC) \to
\oplus_{(C)} KO_n(C^*_rC))
$$
which maps to $\beta$.  In \cite{BGS}, Botvinnik-Gilkey-Stolz show a strong
version of the GLR conjecture for finite cyclic groups, namely that
$ko^+_n(BC) = \ker(A \circ p : ko_n(BC) \to KO_n(C^*_rC))$.   Thus
$\gamma \in \oplus_{(C)} ko_n^+(BC)$, and
hence its image $\beta$ lies in $ko^+_n(B\G)$.  Thus the strong version of the
GLR conjecture for $\G$ is true, namely that
$$ko^+_n(B\G) = \ker(A \circ p : ko_n(B\G) \to KO_n(C^*_r\G)).$$
Hence, by Theorem \ref{kobordism},
our main theorem is proved.
\end{proof}

\begin{remark} {\rm In fact, \cite{BGS} shows that $ko^+_n(BC)$ is
generated by lens spaces and simply-connected spaces, and our  proof
shows that same is true for
$ko^+_n(B\G)$.}
\end{remark}

\begin{remark}
{\rm Fuchsian groups have almost periodic cohomology, i.e. $H^m(\G)\cong
H^{m+2}(\G)$ provided $m>2$.  Thus our main
result is not too surprising in view of \cite{BGS}; the $K$-theory
terms arising from $X_g$ in a sense measures the failure
of $H^*(\G)$ to be periodic. }
\end{remark}

\noindent
{\bf Acknowledgments.}  We thank Diarmuid Crowley and Mark Johnson for
valuable discussions, and Stephan Stolz for suggesting
that we go ahead and do the case of even torsion.  The first author
acknowledges support of the NSF.  The second author acknowledges the support of
an NSF-AWM Travel Grant and
thanks the mathematics departments of The University of Chicago and
Indiana University for their hospitality while work on this
paper was being done.

\end{document}